\documentclass[english,12pt,leqno]{article}
\usepackage{tikz}
\usetikzlibrary{decorations.markings}
\usetikzlibrary{calc}
\usepackage{wrapfig}

\usepackage[T2A]{fontenc}

\usepackage[latin1, russian, utf8]{inputenc}

\usepackage{amsxtra}
\usepackage{amsmath}
\usepackage{amssymb}
\usepackage{amsfonts}
\usepackage[all]{xy}
\usepackage{mathrsfs}
\usepackage{amsthm}
\usepackage{float}
 \usepackage{caption}

\newtheorem{thm}[equation]{Theorem}

\newtheorem{cor}[equation]{Corollary}

\newtheorem{prop}[equation]{Proposition}

\newtheoremstyle{example}{\topsep}{\topsep}%
     {}%         Body font
     {}%         Indent amount (empty = no indent, \parindent = para indent)
     {\bfseries}% Thm head font
     {.}%        Punctuation after thm head
     {2pt}%     Space after thm head (\newline = linebreak)
     {\thmname{#1}\thmnumber{ #2}\thmnote{ #3}}%         Thm head spec

   \theoremstyle{example}
   
   \newtheorem{Defi}[equation]{Definition}
   
   \newtheorem{rem}[equation]{Remark}

   \newtheorem{exas}[equation]{Examples}
   \newtheorem{ex}[equation]{Example}

\newtheoremstyle{example}{\topsep}{\topsep}%
     {}%         Body font
     {}%         Indent amount (empty = no indent, \parindent = para indent)
     {\bfseries}% Thm head font
     {.}%        Punctuation after thm head
     {2pt}%     Space after thm head (\newline = linebreak)
     {\thmname{#1}\thmnumber{ #2}\thmnote{ #3}}%         Thm head spec

   \numberwithin{equation}{section}

\def\CC{\mathbb{C}}
\def\DD{\mathbb{D}}

\def\RR{\mathbb{R}}
\def\ZZ{\mathbb{Z}}

\def\HH{\mathbb{H}}

\def\Aen{\mathfrak{A}}

\def\Pen{\mathfrak{P}}
\def\Qen{\mathfrak{Q}}

\def\Ac{\mathcal{A}}

\def\Cc{\mathcal{C}}

\def\Dc{\mathcal{D}}
\def\Ec{\mathcal{E}}
\def\Fc{\mathcal{F}}

\def\Mc{\mathcal{M}}

\def\Oc{\mathcal{O}}

\def\Rc{\mathcal{R}}

\def\kb{\mathbf{k}}

\def\<{\langle}
\def\>{\rangle}

\def\Act{{\Ac ct}}

\def\Aut{\on{Aut}}

\def\be{\begin{equation}}
\def\ee{\end{equation}}

\def\bef{\begin{figure}[H]\centering}
\def\enf{\end{figure}}

\def\btp{\begin{tikzpicture}}
\def\etp{\end{tikzpicture}}

\def\codim{{\on{codim}}}

\def\Ed{\on{{Ed}}}

\def\Fun{\on{{Fun}}}

\def\Hom{\on{Hom}}

\def\Id{{\on{Id}}}

\def\k {\mathbf k}
\def\Ker{\on{Ker}}

\def\lra{\longrightarrow}
\def\iso{\overset{\sim}{\longrightarrow}}

\def\lra{\longrightarrow}

\def\Ob{\on{{Ob}}}
\def\ol{\overline}
\def\on{\operatorname}

\def\orr{\on{or}}

\def\Perv{\on{Perv}}
\def\PolPerv{\on{PolPerv}}
\def\pt{{\on{pt}}}

\def\Rep{\on{Rep}}

\def\Sh{\on{Sh}}

\def\Surf{\on{Surf}}

\def\ul{\underline}

\def\Vect{\on{Vect}}
\def\Vert{\on{Vert}}

\title{ Perverse sheaves and graphs on surfaces}

\author{ Mikhail Kapranov, Vadim Schechtman}
\begin{document}
%\begin{abstract}

%\end{abstract}

 \maketitle

% \addtocounter{section}{-1}
 %\tableofcontents

\section{Introduction}\label{sec:intro}

The aim of this note is to propose a combinatorial description of the categories  
$\Perv(S,N)$ of perverse sheaves  on a Riemann surface $S$  with possible  singularities 
at a finite set of points $N$. Here we consider 
the sheaves with values in the category   of vector spaces over a fixed base field $\kb$. 
We allow (in fact, require)  the surface $S$ to have a boundary.  Categories of the type $\Perv(S,N)$ can be
used as inductive building blocks in the study of any category of perverse sheaves, see
\cite{gelfand-MV}. Therefore  explicit combinatorial descriptions of them are desirable. 

\vskip .2cm

We proceed in a manner similar to \cite{KS}. 
To define a combinatorial data we need to  fix a {\it Lagrangian skeleton} of $S$. In our case it will be a {\it spanning graph $K\subset S$} 
with the set of vertices $\Vert(K) = N$ (for the precise meaning of the word ''spanning'' see Section \ref{sec:graph-descr}B). (In the case of a hyperplane arrangement in $\CC^n$ discussed in \cite{KS} the Lagrangian skeleton was  $\RR^n\subset \CC^n$.) We denote by $\Ed(K)$ the set of edges of $K$. 
We suppose for simplicity in this Introduction that $K$ has no 
loops. As any graph embedded into an oriented surface, $K$ is naturally a {\em ribbon graph}, i.e., it is equipped with
 a cyclic order on the
set of edges incident to any vertex. 

\vskip .2cm

To any ribbon graph $K$ we associate a category $\Ac_K$  whose 
objects are collections  $\{E_x, E_e\in \Vect(\kb),\ x\in \Vert(K), e\in \Ed(K)\}$ together with linear maps
\[
  \xymatrix{
 E_x  \ar@<.4ex>[r]^{  \gamma}& E_e \ar@<.4ex>[l]^{  \delta}. 
 } 
 \]
given 
for each couple $(x, e)$ with  $x$ being a vertex of an edge $e$.  
Here $\Vect(\kb)$ denotes the category of finite  dimensional $\kb$-vector spaces. 
These 
maps must satisfy the relations which use the ribbon structure on $K$ and are  listed in Section \ref{sec:graph-descr}C  below,.

\vskip .2cm

If $K\subset S$ is a spanning graph as above, then 
our main result (see Theorem \ref{thm:main}) establishes an equivalence 
of categories 
$$
Q_K:\ \Perv(S,N) \iso \Ac_K 
$$
For $\Fc\in \Perv(S,N)$ the vector spaces $Q_K(\Fc)_x, Q_K(\Fc)_v$ are the stalks of the constructible complex
 $\Rc_K(\Fc)=Ri_{K}^!(\Fc)[1]$ on $K$
which, as we prove, is identified with  a constructible sheaf in degree 0. 

\vskip .2cm

A crucial particular case is $S = $ the unit disc  $ D \subset \CC$, 
$N = \{0\}$. Take for the skeleton a corolla $K_n$ with center  at $0$ 
and $n$ branches.  Thus, the same category $\Perv(D,0)$ has infinitely many 
incarnations, being equivalent to $\Ac_n := \Ac_{K_n},\ n\geq 1$. 

The corresponding equivalence $Q_n$ is described 
in Section \ref{sec:fracspin}, see Theorem  \ref{thm:An}.  
The  special cases of this equivalence are:
\begin{itemize}

\item[(i)] $Q_1: \Perv(D,0)\iso \Ac_1$: this is a classical theorem, in the form given in \cite{GGM}. 

\item[(ii)]  $Q_2: \Perv(D,0)\iso \Ac_2$ is a particular case of the main result of  \cite{KS}, see {\it op. cit.},  \S 9A. In {\it loc. cit.} we have also described 
the resulting equivalence
$$
\Ac_1 \iso \Ac_2
$$
explicitely. In a way,  objects of $\Ac_2$ are ''square roots'' of 
objects  of $\Ac_1$, in the same manner as the Dirac operator is 
a square root of the Schr\"odinger operator. That is why we call 
$\Ac_n$ a ''$1/n$-spin (parafermionic) incarnation'' of $\Perv(D,0)$. 
\end{itemize}

\noindent Finally, in Section \ref{sec:graph-descr}D we give (as an easy corollary of the previous discussion) a combinatorial description of 
the category $\PolPerv(S,N)$ of {\it polarized} perverse sheaves, cf. \cite{saito}.  These objects arise ''in nature'' as decategorified
perverse Schobers, cf. \cite{KS-schobers}, the polarization being induced by the 
Euler form  $(X, Y) \mapsto \chi(R\Hom(X,Y))$.

\vskip .2cm
The idea of localization on a Lagrangian skeleton was proposed by M. Kontsevich in the context of Fukaya categories.
The fact that it is also applicable to the problem of classifying perverse sheaves
(the constructions of \cite{GGM} and \cite{KS} can be seen, in retrospect, as manifestations of this idea)
is a remarkable phenomenon.  It
indicates a deep connection between Fukaya categories and perversity. A similar approach will be used in
\cite{DKSS} to construct the Fukaya category of a surface with coefficients in a perverse Schober.

The work of M.K.  was supported by World Premier International Research Center Initiative (WPI Initiative), MEXT, Japan.

%\vfill\eject

 \section{The ``fractional spin'' description of perverse sheaves on the disk}\label{sec:fracspin}
 
 \noindent {\bf A. Statement of the result.} 
 Let $X$ be a complex manifold. By a {\em perverse sheaf} on $X$ we mean a $\CC$-constructible complex $\Fc$ of
 sheaves of $\k$-vector spaces on $X$, 
 which satisfies the middle perversity condition, normalized so that a local system in degree $0$ is perverse.
 Thus, if  $\k=\CC$ and $\Mc$ is a holonomic $\Dc_X$-module,   then $\ul{R\Hom}_{\Dc_X}(\Mc, \Oc_X)$ is a perverse sheaf.
 
  \vskip .2cm

Let $D=\{|z| < 1\} \subset \CC$ be the unit disk. Let $\k$ be a field. 
 We denote by $\Perv(D,0)$ the abelian category of perverse shaves of $\k$-vector
 spaces on $D$ which are smooth (i.e., reduce to a local system in degree $0$)
 outside $0\in D$.

 Let $n\geq 1$ be an integer.  Let  $\Ac_n$ be the category of diagrams
 of finite-dimensional $\k$-vector spaces (quivers) $Q$ , consisting of   spaces
 $E_0, E_1, \cdots, E_n$ and linear maps
 \[
 \xymatrix{
 E_0 \ar@<.4ex>[r]^{\gamma_i}& E_i \ar@<.4ex>[l]^{\delta_i}
 }, \quad i=1, \cdots, n, 
 \]
  satisfying the conditions (for $n\geq 2$): 
  \begin{itemize}
  \item[(C1)] $\gamma_i\delta_i = \Id_{E_i}$. 
  
  \item[(C2)] The operator 
  $T_i := \gamma_{i+1}\delta_i: E_i \to E_{i+1}$ (where $i+1$ is considered modulo $n$), is an
  isomorphism for each $i=1,\cdots, n$. 
  
 \item[(C3)]  For $i\neq j, j+1 \mod n$, we have $\gamma_i\delta_j=0$. 
  
  \end{itemize}

  For $n=1$ we impose the standard relation:
  
  \begin{itemize}
  \item[(C)] The operator $T=\Id_{E_1} -\gamma_1\delta_1:  E_1\to E_1$ is an isomorphism. 
  \end{itemize}

  \begin{thm}\label{thm:An}
  For each $n\geq 1$, the category $\Perv(D,0)$ is equivalent to $\Ac_n$. 
  
  \end{thm}
  
  For $n=1$ this is the standard $(\Phi, \Psi)$ description of perverse sheaves on the disk
  \cite{GGM}, \cite{beil-gluing}. For $n=2$ this is a particular case of the description in
  \cite{KS} (\S 9 there). 
  
  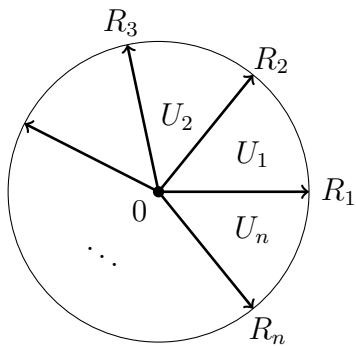
\begin{wrapfigure}{r}{6cm}
 \centering
    \begin{tikzpicture}[scale=0.5]
    
    \node (0) at (0,0){};
    \fill (0) circle (0.15); 
    
    \draw (0) circle (4.0); 
    
    \draw[->, line width = 1] (0,0) -- (4,0);
    
    \node (1) at (4,0){};
  %  \fill (1) circle (0.15); 

    \node at (-0.5, -0.5) {$0$}; 
    \node at (4.5, -0.5) {$$}; 
    
    \draw [->, line width = 1](0,0) -- (51:4); 
    
    \draw [->, line width = 1](0,0) -- (102:4); 
      
    \draw [->, line width = 1](0,0) -- (153:4); 
         
   \draw [->, line width = 1](0,0) -- (-51:4); 
   
   \node at (-1.5,-1.5){$\ddots$};

   \node at (4.8,0){$R_1$}; 
   
   \node at (3,3.5 ){$R_2$}; 
   
   \node at (-1,4.5){$R_3$}; 
   
   \node at (2.9, -3.7){$R_n$}; 
   
   \node at (2.5,-1){$U_n$}; 
   
   \node at (2.5,1){$U_1$}; 
   
   \node at (0.5,2){$U_2$}; 
   
    \end{tikzpicture}
   \caption{The graph $K=K_n$. }\label{fig:corolla}
 \end{wrapfigure}

 \noindent {\bf B. Method of the proof.}  For the proof we consider a star shaped graph $K=K_n\subset D$ obtained by drawing $n$ radii $R_1, ..., R_n$
from $0$ in the counterclockwise order, see Fig. \ref{fig:corolla}. Then $D-K$ is the union of $n$ open sectors $U_1, \cdots, U_n$ numbered
so that $U_\nu$ is bordered by $R_\nu$ and $R_{\nu+1}$. 

  \begin{prop}\label{prop:pure}
  For any $\Fc\in\Perv(D,0)$ we have $\ul\HH^i_K(\Fc)=0$ for $i\neq 1$. Therefore the functor
  \[
  \Rc: \Perv(D,0) \lra \Sh_K, \quad \Fc\mapsto \Rc(\Fc) = \Rc_K(\Fc):= \ul\HH^1_K(\Fc)
  \]
  is an exact functor of abelian categories. 
  \end{prop}

 \noindent {\sl Proof:} Near a point $x\in K$ other than $0$, the graph $K$ is
 a real codimension 1 submanifold in $D$, and $\Fc$ is a local system in degree 0,
 so the statement is obvious (``local Poincar\'e duality''). So we really need only
  to prove that the space
 $\ul\HH^j_{K_n}(\Fc)_0 = \HH^{j}_{K_n}(D,\Fc)$  vanishes for $j\neq 1$. 
The case $n=0$ is known, $\ul\HH^1_{K_1}(D,\Fc)$  being
identified with $\Phi(\Fc)$, see \cite{GGM}. The general case is proved by
induction on $n$. We consider an embedding  of $K_n$ into $K_{n+1}$
so that the new radius $R_{n+1}$ subdivides the sector $U_n$ into two.
 This leads to  a morphism between the long exact sequences  
 relating hypercohomology with and without support in  $K_n$ and $K_{n+1}$:
\[
\xymatrix{
\cdots   \ar[r]&  \HH^j(D-K_n, \Fc) \ar[d]_{\alpha_j} \ar[r] &\HH^{j+1}_{K_n}(D,\Fc)
\ar[d]^{\beta_j}\ar[r]& \HH^{j+1}(D,\Fc) \ar[d]^{=}\ar[r] &\cdots
\\
\cdots  \ar[r]&  \HH^j(D-K_{n+1}, \Fc)  \ar[r] &\HH^{j+1}_{K_{n+1}}(D,\Fc)
 \ar[r]& \HH^{j+1}(D,\Fc) \ar[r] &\cdots
}
\]
For $j\neq 0$ the map $\alpha_j$ is an isomorphism because its source and target are $0$. Indeed, $D-K_n$ as
well as $D-K_{n+1}$ is the union of contractible sectors, and $\Fc$ is a local system in degree $0$ outside $0$, so the
higher cohomology of each sector with coefficients in $\Fc$ vanishes. This means that 
\[
\HH^{j+1}_{K_n+1}(D, \Fc) \,\,=\,\, H^{j+1}_{K_n}(D, \Fc) \,\,=\,\,0, \quad j\neq 0
\]
and so by induction all these spaces are equal to 0.
 \qed
 
 \begin{rem}
 An alternative proof of  the vanishing of $\ul\HH^{\neq 1}_{K_n}( \Fc)_0$ can be obtained by noticing that
 $\ul{R\Gamma}_{K_n}(\Fc)_0[1]$ can be identified with $\Phi_{z^n}(\Fc)$, the space of vanishing cycles
 with respect to the function $z^n$. It is known that forming the sheaf of vanishing cycles with respect to
 any holomorphic function preserves perversity. 
 \end{rem} 
 
 The graph $K=K_n$ is a regular cellular space with cells  being $\{0\}$ and the open rays $R_1, \cdots, R_n$. 
  For $\Fc\in\Perv(D,0)$ the sheaf $\Rc(\Fc)$ is a cellular sheaf on $K$ and as such is completely determined by
  the linear algebra data of:
  
  \begin{itemize}
  
  \item[(1)] Stalks at the (generic point of the) cells, which we denote:
  \[
  \begin{gathered}
  E_0 = E_0(\Fc_n) := \Rc(\Fc)_0 = \text{ stalk at } 0; \\
  E_i = E_i(\Fc) = \Rc(\Fc)_{R_i} = \text{ stalk at } R_i, \quad i=1, \cdots, n. 
  \end{gathered}
  \]
  
  \item[(2)] Generalization maps corresponding to inclusions of closures of the cells, which we denote
  \[
  \gamma_i: E_0 \lra E_i, \quad i=1, \cdots, n. 
  \]
  
  \end{itemize}
  
  This gives ``one half'' of the quiver we want to associate to $\Fc$. 
  
  \vskip .2cm
  
  \noindent {\bf C. Cousin complex.} In order to get the second half of the maps (the $\delta_i$), we introduce,
  by analogy with \cite{KS}, a canonical ``Cousin-type'' resolution of any $\Fc\in\Perv(D,0)$.
  
  Denote by
  \[
  i: K\hookrightarrow D, \quad j: D-K\hookrightarrow D
  \]
  the embeddings of the closed subset $K$ and of its complement $D-K=\bigsqcup_{\nu=1}^n U_\nu$. 
  For any complex of sheaves $\Fc$ on $D$ (perverse or not) we have a canonical
  distinguished triangle in $D^b\Sh_D$:
  \be\label{eq:triangle}
  i_* i^!\Fc \lra \Fc\lra j_* j^*\Fc \buildrel \delta \over \lra i_* i^!\Fc[1]
  \ee
  (here and elsewhere $j_*$ means the full derived direct image). Recalling that $i^!$ has the meaning of
  cohomology with support, and denoting $j_\nu: U_\nu\hookrightarrow D$ the embeddings of the connected
  components of $D-K$, we conclude, from Proposition \ref{prop:pure}:

  \begin{cor}
  Let $\Fc\in\Perv(D,0)$. Then $\Fc$ is quasi-isomorphic (i.e., can be thought of as represented by)
  the following 2-term complex of sheaves on $D$:
  \[
  \Ec^\bullet(\Fc) \,\,=\,\,\biggl\{ \bigoplus_{\nu=1}^n j_{\nu *} \bigl(\Fc|_{U_\nu}\bigr) 
  \buildrel \ul\delta\over\lra \Rc(\Fc)\biggr\}.
  \]
  Here the grading of the complex is in degrees $0,1$ and the map $\ul\delta$ is induced by the boundary
  morphism $\delta$ from \eqref{eq:triangle}. \qed
   \end{cor}
  
  We further identify
   \[
  \Fc|_{U_\nu} \simeq{ {\ul E}_\nu}_{U_\nu} = \text{ constant sheaf on } U_\nu \text{ with stalk } E_\nu
  \]
  as follows. As $\Fc$ is locally constant (hence constant) on $U_\nu$, it is enough to specify an isomorphism
  \[
  H^0(U_\nu, \Fc) \buildrel \alpha_\nu\over\lra E_\nu = \ul\HH^1_K(\Fc)_x,
  \]
  where $x$ is any point on $R_\nu$. 
  Taking a small disk $V$ around $x$, we have $E_\nu= \HH^1_{V\cap K}(V,\Fc)$. This hypercohomology with support 
   is identified, via the coboundary map of the standard long exact sequence relating cohomology with and without support),
    with the 
 cokernel of the map
 \[
 H^0(V, \Fc) \lra H^0(V\cap U_\nu, \Fc) \oplus H^0(V\cap U_{\nu-1}, \Fc) =  H^0(V-K, \Fc)
 \]
 Since $\Fc$ is constant on $V$, the projection of $H^0(V\cap U_\nu, \Fc)=H^0(U_\nu, \Fc)$ to the cokernel, i.e., to $E_\nu$,
 is an isomorphism. We define $\alpha_\nu$ to be this projection. 
  
  \vskip .2cm
  
  We now assume $n\geq 2$. Then the closure of $U_\nu$ is a proper closed sector $\ol U_\nu$,
  and so we can rewrite the complex $\Ec^\bullet(\Fc)$ as 
  \[
  \Ec^\bullet(\Fc) \,\,=\,\,\biggl\{ \bigoplus_{\nu=1}^n  \ul {E_\nu}_{\, \ol U_\nu}
  \buildrel \ul\delta\over\lra \Rc(\Fc)\biggr\}.
  \]
 Indeed,
 \[
 j_{\nu *}(\Fc|_{U_\nu} = j_{\nu *}({\ul E_\nu}_{U_\nu} ) \,\,=\,\,  \ul {E_\nu}_{\, \ol U_\nu}. 
 \]

  \vskip .2cm
  
  \noindent {\bf D. Analysis of the morphism $\ul\delta$.}
We now analyze the maps of stalks over various points induced by $\ul\delta$. Since the target of $\ul\delta$ is supported on $K$, it is
enough to consider two cases:

\vskip .2cm

\ul{Stalks over $0$.} We get maps of vector spaces
\[
\delta_\nu: \bigl({\ul E_\nu}_{\ol U_\nu}\bigr)_0 = E_\nu \lra \Rc(\Fc)_0 = E_0, \quad \nu=1, \cdots, n. 
\]
These maps, together with the generalization maps $\gamma_\nu$, form the quiver
\[
Q=Q(\Fc)  \,\,=\,\, \bigl\{ 
\xymatrix{
 E_0 \ar@<.4ex>[r]^{\gamma_i}& E_i \ar@<.4ex>[l]^{\delta_i}
 }, \quad i=1, \cdots, n 
 \bigr\}
\]
which we associate to $\Fc$.

\vskip .2cm

\ul{Stalks over a generic point $x\in R_\nu$.} As $x$ lies in two closed subsets $\ul U_\nu$ and $\ul U_{\nu -1}$, we have two maps
\[
\begin{gathered}
\delta_{U_\nu, R_\nu}: \bigl({\ul E_\nu}_{\ol U_{\nu}}\bigr)_x = E_{\nu} \lra \Rc(\Fc)_x = E_\nu, 
\\
\delta_{U_{\nu-1}, R_\nu}: \bigl({\ul E_{\nu-1}}_{\ol U_{\nu-1}}\bigr)_x = E_{\nu-1} \lra \Rc(\Fc)_x = E_\nu. 
\end{gathered}
\] 
  
  \begin{prop}\label{prop:rels}
 We have the following relations:
 \[
 \delta_{U_\nu, R_\nu} = \Id_{E_\nu}, \quad \delta_{U_{\nu -1}, R_\nu}=-T_{\nu -1}, 
 \]
 where 
 \[
 T_{\nu -1}: E_{\nu -1} = h^0(U_{\nu -1}, \Fc) \lra E_\nu = H^0(U_\nu, \Fc)
 \]
 is the counterclockwise monodromy map for the local system $\Fc|_{D-\{0\}}$.   
  \end{prop}
  
 \noindent {\sl Proof:} We start with the first relation. Recal that the identification $\alpha_\nu: H^0(U_\nu, \Fc)\to E_\nu$
 was defined in terms of representation of $E_\nu$ as a quotient, i.e., in terms of the coboundary map in the
  LES relating hypercohomology with and without support. So the differential $\ul\delta$ in $\Ec^\bullet(\Fc)$,
  applied to a section $s\in H^0(U_\nu, \Fc)$, gives precisely $\alpha_\nu(s)$, so after the identification by $\alpha_\nu$,
  the map on stalks over $x\in R_\nu$, becomes the identity.
  
  \vskip .2cm
  
  We now prove the second relation. Representing $e\in E_{\nu -1}$ by a section $s\in H^0(U_{\nu -1}, \Fc)$,
  we see that $\delta_{U_{\nu -1}, R_\nu}(e)$ is represented by the image of
  \[
  (0,s) \,\,\in \,\, H^0(U_\nu, \Fc) \oplus H^0(U_{\nu -1}, \Fc)
  \]
  in the quotient
  \[
  \bigl( H^0(U_\nu, \Fc) \oplus H^0(U_{\nu -1}, \Fc)\bigr) \bigl / H^0(V, \Fc).
  \]
  But the identification of this quotient with $E_\nu$ is via the projection to the first, not second, summand, i.e., to
  $H^0(V\cap U_\nu, \Fc)=H^0(U_\nu, \Fc)$. The element $(t,0)$ projecting to the same element of the quotient as $(0,s)$,
  has $t=-T_{\nu -1}(s)=$ {\em minus} the analytic continuation of $s$ to $U_\nu$.\qed

  \begin{prop}
  The maps $\gamma_\nu, \delta_\nu$ in the diagram $Q(\Fc)$ satisfy the conditions (C1)-(C3), i.e., $Q(\Fc)$
  is an object of the category $\Ac_n$. 
  \end{prop}
  
  \noindent {\sl Proof:|} We spell out the conditions that the differential $\ul\delta$ in the Cousin complex
  $\Ec^\bullet(\Fc)$ is a morphism of sheaves. More precisely, both terms of the complex are cellular sheaves on $D$
  with respect to the regular cell decomposition given by $0, R_1,\cdots, R_n, U_1, \cdots, U_n$. 
  So the maps of the stalks induced by $\ul\delta$, must commute with the generalization maps.
  
  Consider the generalization maps from $0$ to $R_\nu$. In thje following diagram the top row is the stalk of the
  complex $\Ec^\bullet(\Fc)$ over $0$, the bottom row is the stalk over $R_\nu$, and the vertical arrows are the
  generalization maps:
  \[
  \xymatrix{
  \bigoplus_ {\mu=1}^n E_\mu 
  \ar[d]_{p_{\nu, \nu-1}}
   \ar[rr]^{\sum \delta_\mu}  && E_0
   \ar[d]^{\gamma_\nu}
  \\
  E_\nu \oplus E_{\nu -1} \ar[rr]_{\hskip .3cm \Id - T_{\nu-1}}&& E_\nu,
  }
  \]
  the lower horizontal arrow having been described in Proposition \ref{prop:rels}. We now spell out the condition of
  commutativity on each summand $E_\mu$ inside $\bigoplus_{\mu=1}^n E_\mu$. 
  
  \vskip .2cm
  
  \ul{Commutativity on $E_\nu$:} this means $\gamma_\nu\delta_\nu=\Id$.
  
  \vskip .2cm
  
  \ul{Commutativity on $E_{\nu -1}$:} this means $\gamma_\nu \delta_{\nu -1}=-T_{\nu-1}$, in particular, this composition
  is an isomorphism. 
  
  \vskip .2cm
  
  \ul{Commutativity on $E_\mu$ for $\mu  \neq \nu, \nu -1$:} This means that $\gamma_\nu\delta_\mu=0$, since the
  projection $p_{\nu, \nu -1}$ annihilates $E_\mu$. The proposition is proved. \qed
  
  \vskip .3cm
  
  \noindent {\bf E. $Q(\Fc)$ and duality.} Recall that the category $\Perv(D,0)$ has  a perfect duality
  \be\label{eq:star}
  \Fc \mapsto \Fc^\bigstar = \DD(\Fc)[2],
  \ee
  i.e., the shifted Verdier duality normalized so that for $\Fc$ being a  local system (in our case, constant sheaf) in degree $0$, we have that
  $\Fc^\bigstar$ is the dual local system in degree 0. We will use the notation \eqref{eq:star}  also for
  more general complexes of sheaves on $D$. 
  
  On the other hand, the category $\Ac_n$ also has a perfect duality
  \[
  Q =  \bigl\{ 
\xymatrix{
 E_0 \ar@<.4ex>[r]^{\gamma_\nu}& E_\nu\ar@<.4ex>[l]^{\delta_\nu}
 }  \bigr\}_{\nu=1}^n \,\,\, \mapsto \,\,\, Q^* =  \bigl\{ 
\xymatrix{
 E_0^* \ar@<.4ex>[r]^{\delta^*_\nu}& E_\nu^* \ar@<.4ex>[l]^{\gamma^*_\nu}
 }  \bigr\}_{\nu=1}^n.
  \]
  
  \begin{prop}\label{prop:verdier}
  The functor $Q: \Perv(D,0)\to\Ac_n$ commutes with duality, i.e., we have canonical identifications
  $Q(\Fc^\bigstar) \simeq Q(\Fc)^*$. 
    \end{prop}
    
    \noindent {\sl Proof:} We modify the argument of \cite{KS}, Prop. 4.6.
    That is, we think of $K$ as consisting of $n$ ``equidistant'' rays $R_\nu$, joining $0$ with
    \[
    \zeta^{\nu-1}, \,\,\, \zeta = e^{2\pi i/n}, \,\,\, \nu = 1,2, \cdots, n. 
    \]
    We consider another star-shaped graph $K'$ formed by the radii $R'_1, \cdots, R'_n$ so that
    $R'_\nu$ is in the middle of the sector $U_\nu$. Thus, the rotation by $e^{i\pi/n}$ identifies $R_\nu$ with $R'_\nu$\
    and $K'$ with $K$. 
    
    We can use $K'$ instead of $K$ to define $\Rc(\Fc)$ and $\Rc(\Fc^\bigstar)$.We will denote the corresponding sheaves
    $\Rc_{K'}(\Fc) = \ul\HH^1_{K'}(\Fc)$, and similarly for $\Rc_{K'}(\Fc^\bigstar)$. The sheaves defined using $K$, will
    be denoted by $\Rc_K(\Fc)$ and so on. 
    
    Since Verdier duality interchanges $i^!$ and
    $i^*$ (for $i: K'\to D$ being the embedding), we have
    \[
    \Rc_{K'}(\Fc^\bigstar)^\bigstar \,\,\simeq \,\, \Fc|_{K'}
    \]
    (usual restriction). To calculate this restriction, we use the Cousin resolution of $\Fc$ defined by using $K$ and the $U_\nu$:
    \[
    \Fc \,\,\simeq \,\, \Ec^\bullet \,\,=\,\,
    \biggl\{ \bigoplus_{\nu=1}^n {\ul E_\nu}_{\ol U_\nu}\buildrel\ul\delta\over\lra \Rc_K(\Fc)\biggr\}.
    \]
  So we restrict $\Ec^\bullet$ to $K'$. SInce $K'\cap K=\{0\}$ and $\Rc_K(\Fc)$ is supported on $K$, the restriction
  $\Rc_{K}(\Fc)|_{K'} = {\ul E_0}_0$ is the skyscraper sheaf at $0$ with stalk $E_0$. So
  \[
  \Fc|_{K'} \,\,\sim \,\, \biggl\{ \bigoplus_{\nu=1}^n ({\ul E_\nu})_{\ol R'_\nu} \buildrel \ul\delta'|_{K'} = \sum\delta_\nu\over\lra ({\ul E_0})_0\biggr\}. 
  \]
  On the other hand, the shifted Verdier dual to $\Rc_{K'}(\Fc^\bigstar)$, as a sheaf on $K'$ is identified by, e.g., \cite{KS}, Prop. 1.11
   with the complex of sheaves
  \[
 \biggl\{  \bigoplus_{C\subset K' \atop \codim (C)=0} \ul{E_C(\Fc^\bigstar)^*}_{\ol C} \otimes \orr(C) 
  \buildrel \sum  (\gamma_{C'C}^{\Fc^\bigstar})^*\over\lra  \bigoplus_{C\subset K' \atop \codim (C)=1} \ul{E_C(\Fc^\bigstar)^*}_{\ol C} \otimes \orr(C) 
  \biggr\}. 
  \]
  Here $C$ runs over all cells of the cell complex $K'$, and $E_C(\Fc^\bigstar)$ is the stalk of the cellular sheaf 
  $\Rc_{K'}(\Fc^\bigstar)$ at the cell $C$.
  Explicitly, $C$ is either $0$ or one of the $R'_\nu$, so 
  \[
  \Rc_{K'}(\Fc^\bigstar)^\bigstar \,\,=\,\,\biggl\{
  \bigoplus_{\nu=1}^n \ul{E_\nu(\Fc^\bigstar)^*}_{\ol R_\nu}
   \buildrel \sum  (\gamma_{C'C}^{\Fc^\bigstar})^*\over\lra  
  \ul{E_0(\Fc^\bigstar)^*}_0\biggr\}. 
  \]
  By the above, this complex is quasi-isomorphic to
  \[
 \Fc |_{K'} \,\,= \,\, \biggl\{ \bigoplus_{\nu=1}^n \ul {E_\nu(\Fc)}_{\ol {R_\nu}} \buildrel \sum \delta_\nu^\Fc \over\lra
  \ul{E_0(\Fc)}_0\biggr\}.
  \]
  So we conclude that
  \[
  E_\nu(\Fc^\bigstar) \,\,=\,\, E_\nu(\Fc)^*, \quad \gamma_\nu^{\Fc^\bigstar} =(\delta_\nu^\Fc)^*. 
  \]
  This proves the proposition.
  
   \vskip .3cm
  
  \noindent {\bf Proof of Theorem \ref{thm:An}. } We already have the functor
  \[
  Q: \Perv(D,0)\lra\Ac_n, \quad \Fc\mapsto Q(\Fc). 
  \]
Let us define a functor $\Ec: \Ac_n\to D^b\Sh_D$. Suppose we are given   
\[
Q =  \bigl\{ 
\xymatrix{
 E_0 \ar@<.4ex>[r]^{\gamma_\nu}& E_\nu\ar@<.4ex>[l]^{\delta_\nu}
 }  \bigr\}_{\nu=1}^n \in\Ac_n.
 \]
 We associate to it the Cousin complex
 \[
 \Ec^\bullet(Q) \,\,=\,\,\biggl\{ \bigoplus_{\nu=1}^n \ul{E_\nu}_{\ol U_\nu} 
 \buildrel \ul\delta\over\lra \Rc(Q)\biggr\}.
 \]
 Here $\Rc(Q)$ is the cellular sheaf on $K$ with stalk $E_0$ at $0$, stalk $E_\nu$ at $R_\nu$
 and the generalization map from $0$ to $R_\nu$ given by $\gamma_\nu$. The map $ul\delta$ 
  is defined on the stalks as follows:
  
  \vskip .2cm
  
  \ul{Over $0$:} the map 
  \[
  (\ul {E_\nu}_{\ol U_\nu})_0 = E_\nu  \lra \Rc(Q)_0=E_0
  \]
   is given by $\delta_\nu: E_\nu\to E_0$.

  \vskip .2cm
  
  \ul{Over $R_\nu$:} The map 
  \[
  \biggl( \bigoplus_{\mu=1}^n \ul{E_\mu}_{\ol U_\mu}\biggr)_{R_\nu} = \,\, E_\nu \oplus E_{\nu -1} 
  \lra\,\,  \Rc(Q)_{R_|nu}=E_\nu
  \]
  is given by
  \[
  \Id - T_{\nu -1}: E_\nu \oplus E_{\nu -1} \lra E_\nu. 
  \]
  Reading the proof of Proposition \ref{prop:rels} backwards, we see that the conditions (C1)-(C3) mean that
  in this way we get a morphism $\ul\delta$ of cellular sheaves on $D$, so $\Ec^\bullet(Q)$ is an
  object of $D^b\Sh_D$. 
  
  Further, similarly to Proposition \ref{prop:verdier}, we see that $\Ec^\bullet(Q^*) \simeq \Ec(Q)^\bigstar$.
  
  \begin{prop}\label{prop:conper}
  $\Ec^\bullet(Q)$ is constructible with respect to the stratification $(\{0\}, D-\{0\})$ and is perverse. 
  \end{prop}
  
  \noindent {\sl Proof:} \ul{Constructibility.} It is sufficient to prove the following:
  \begin{itemize}
  \item[(a)] The sheaf $\ul H^0(\Ec^\bullet(Q))|_{D-\{0\}}$ is locally constant. 
  
  \item[(b)] The sheaf  $\ul H^1(\Ec^\bullet(Q))|_{D-\{0\}}$ is equal to $0$. 
  \end{itemize}
  
  To see (a), we look at the map of stalks over $R_\nu$:
  \[
  \Id-T_{\nu -1}: E_\nu\oplus E_{\nu -1} \lra E_\nu
  \]
  given by the differential $\ul\delta$ in $\Ec^\bullet(Q)$. So, by definition, $\ul H^0(\Ec^\bullet(Q))_{R_\nu}$
  and $\ul H^1(\Ec^\bullet(Q))_{R_\nu}$ are the kernel and cockerel of this map. 
  
  Now, since $T_{\nu -1}$ is an isomorphism, $\Ker (\Id - T_{\nu -1})$ projects to both $E_\nu$ and $E_{\nu -1}$
  isomorphically. This means that $\ul H^0(\Ec^\bullet(Q))$ is locally constant over $R_\nu$:
  the stalk at $R_\nu$ projects (``generalizes'') to the stalks at $U_\nu$ and $U_{|nu -1}$ in an isomorphic way.
  
  To see (b), we notice that $\Id-T_{\nu -1}$ is clearly surjective and so 
  $\ul H^1(\Ec^\bullet(Q))_{R_\nu}=0$. Since $\Ec^1(Q)=\Rc(Q)$ is supported on $K$, this means that 
  $\ul H^1(\Ec^\bullet(Q))$ is supported at $0$, and its restriction to $D-\{0\}$ vanishes. So $\ul H^i(\Ec^\bullet(Q))$
  are $\CC$-constructible as claimed.
  
  \vskip .2cm
  
  \ul{Perversity.}
  By the above, $\Ec^\bullet(Q)$ is {\it semi-perverse}, i.e., lies in the non-positive part of the perverse t-structure, that is,
  $\ul H^i(\Ec^\bullet(Q))$ is supported on complex codimension $\geq i$. 
  Further, $\Ec^\bullet(Q)^\bigstar$ also satisfies the same semi-perversity since it is identified with
  $\Ec^\bullet(Q^*)$. This means that $\Ec^\bullet(Q)$ is fully perverse. Proposition \ref{prop:conper}
  is  proved. 
  
  \vskip .2cm
  
 It remains to show that the functors
  \[
  \xymatrix{
 \Perv(D,0)  \ar@<.4ex>[r]^{\hskip .5cm Q}& \Ac_n\ar@<.4ex>[l]^{\hskip .5cm \Ec}
 } 
  \]
 are quasi-inverse to each other. This is done in a way completely parallel to \cite{KS}, Prop. 6.2 and Lemma
 6.3 (``orthogonality relations''). 
 Theorem  \ref{thm:An} is proved.  $\square$

 \vskip .2cm

See Appendix for some further study of the categories $\Ac_n$. 
  
%\vfill\eject

\section{ The graph description of perverse sheaves on an oriented surface}\label{sec:graph-descr}

 \noindent {\bf A. Generalities. The purity property.} Let $S$ be a compact topological surface, possibly
 with boundary $\partial S$; we denote $S^\circ=S-\partial S$ the interior. 
 Let $N\subset S^\circ$ be a finite subset. We then have the category $\Perv(S,N)$ formed by perverse
 sheaves of $\k$-vector spaces on $S$, smooth outside $N$. 
 
 By a {\em graph} we mean a topological space obtained from a finite 1-dimensional CW-complex
 by removing finitely many points. Thus we do  allow edges
 not terminating in a vertex on some side  ( ``legs''),  as well as 1-valent and 2-valent vertices as well as loops. 
 For a vertex $x$ of a graph $K$ we denote by $H(x)$ the set of
 half-edges incident to $x$. We can, if we wish, consider any point $x\in K$ as a vertex:
 if it lies on an edge, we consider it as a 2-valent vertex, so $H(x)$ is this case is the
 set of the two orientations of the edge containing $x$. 
 Further, for a graph $K$ we denote by $\Vert(K)$ and $\Ed(K)$ the sets of vertices and edges of $K$.

 We denote by $\Cc_K$ the {\em cell category} of $K$ defined as follows.
 The set $\Ob(\Cc_K)$ is $\Vert(K)\sqcup\Ed(K)$ (``cells"). Non-identity morphisms can exist only between
 a vertex $x$ and an edge $e$, and
 \[
 \Hom(x,e) \,\,=\,\, \bigl\{ \text{half-edges } h\in H(x) \text{ contained in } e\bigr\}. 
 \]
 So $|\Hom(x,e)|$ can be $0, 1$ or $2$ (the last posibility happens when $e$ is a loop beginning and ending at $x$). 
 If $K$ has no loops, then $\Cc_K$ is a poset. We denote by $\Rep(\Cc_K)=\Fun(\Cc_K, \Vect_\k)$
 the category of representations of $\Cc_K$ over $\k$. 
 
 \begin{prop}\label{prop:celgraph}
 The category of cellular sheaves on $K$ is equivalent to $\Rep(\Cc_K)$.  
 \end{prop}
 
 \noindent {\sl Proof:} This is a particular case of general statement \cite{treumann} which describes
 constructible sheaves on any stratified space in terms of  representations of the category of exit paths. \qed
 
 \vskip .2cm
 
 Let now  $K\subset S$ be any embedded graph (possibly passing through some points of $N$). 
 We allow 1-valent vertices of $K$ to be  situated inside $S$, as well as on $\partial S$. 
 
 \begin{prop}
 For $\Fc\in\Perv(S,N)$ we have $\ul\HH^i_K(\Fc)=0$ for $i\neq 1$.
 \end{prop}
 
 \noindent {\sl Proof:} This follows from Proposition \ref{prop:pure}. Indeed, the statement is local on $S$, and the graph $K$
 is modeled, near each of its points, by a star shaped graph in a disk. \qed
 
 \vskip .2cm
 
 We denote
 \[
 \Rc(\Fc) \,\,=\,\,\Rc_K(\Fc) \,\, := \,\, \ul\HH^1_K(\Fc),
 \]
 this is a cellular sheaf on $K$.

 \vskip .2cm
 
 \noindent {\bf B. Spanning ribbon graphs.} From now on we assume that $S$ is oriented. Graphs embedded into $S$
 have, therefore, a canonical {\em ribbon structure}, i.e.,  a choice of a cyclic ordering on each set $H(x)$.
 See, e.g., \cite {DK} for more background on this classical concept. 
 
 For a ribbon graph $K$ we have a germ of an oriented surface with boundary $\Surf(K)$ obtained by thickening each edge
 to a ribbon and gluing the ribbons at vertices according to the cyclic order. In the case of a 1-valent vertex $x$ we take the
ribbon  to contain $x$, so that $x$ will be inside $\Surf(K)$. 
 
 By a {\em spanning graph} for $S$ we mean a graph $K$, embedded into $S^\circ$ as a closed subset,
 such that the closure $\ol K\subset S$ is a graph embedded into $S$, and the embedding
  $K\hookrightarrow S^\circ$
 is  a  homotopy equivalence. Thus we allow for  legs of $K$ to touch the boundary of $S$.

 \vskip .2cm
 
 \noindent {\bf C.  The category associated to a ribbon graph.} 
 
 \begin{Defi}
 
 Let $K$ be a  graph, and $\Cc_K$ be its cell category. By a {\em double
 representation} of $\Cc_K$ we mean a datum $Q$ of:
 
 \begin{itemize}
 \item[(1)] For each $x\in\Vert(K)$, a vector space $E_x$.
 
 \item[(2)] For each $e\in\Ed(K)$, a vector space $E_e$.
 
 \item[(2)] For each half-edge $h$ incident to a vertex $x$ and an edge $e$, linear maps
 \[
  \xymatrix{
 E_x  \ar@<.4ex>[r]^{  \gamma_h}& E_e \ar@<.4ex>[l]^{  \delta_h}. 
 } 
 \]
 \end{itemize}
 Let $\Rep^{(2)}\Cc_K$ be the categry of double representations of $\Cc_K$. 
 \end{Defi}
 
  Let now $K$ be a ribbon graph. 
 Denote by $\Ac_K$ the full subcategory in $\Rep^{(2)}\Cc_K$ formed by double representations
 $Q=(E_x, E_e, \gamma_h, \delta_h)$ such that for each vertex $x\in K$ the following conditions are satisfied
 (depending on the valency of $x$):
 
 \begin{itemize} 
 \item If $x$ is 1-valent, then we require:
 
 \begin{itemize}
 \item[($C_x$)] $\Id_{E_e}-\gamma_h\delta_h: E\to E$ is an isomorphism. 
 \end{itemize}
 
 \item If the valency of $x$ is $\geq 2$, then we require:
 
 \begin{enumerate}
 \item[($C1_x$)] For each half-edge $h$ incident to $x$, we have $\gamma_h \delta_h=\Id_{E_e}$.
 
 \item[($C2_x$)] Let $h, h'$ be  any two half-edges incident to $x$ such that $h'$ {\em immediately follows}
 $h$ in the cyclic order on $H(x)$. Let $e, e'$ be the edges containing $h, h'$. 
 Then $\gamma_{h'}\delta_h: E_e\to E_{e'}$ is an isomorphism. 
 
 \item[($C3_x$)] If $h, h'$ are two half edges incident to $x$ such that $h\neq h'$ and $h'$ does not
 immediately follow $h$, then $\gamma_{H'}\delta_h=0$. 
 
 \end{enumerate}

 \end{itemize}
 
 \begin{ex}
 If $K=K_n$ is a ``ribbon corolla''  with one vertex and $n$ legs, then $\Ac_K=\Ac_n$ is the category from \S 
 \ref{sec:fracspin}. 
 \end{ex}

 \vskip .2cm
 
 \noindent {\bf D. Description of $\Perv(S,N)$ in terms of spanning graphs.} Let $S$ be an oriented surface,
 $K\subset S$ be a spanning graph and $N=\Vert(K)$. For $\Fc\in\Perv(S,N)$ we have the sheaf $\Rc_K(\Fc)$
 on $K$, cellular with respect to the cell structure given by the vertices and edges. Therefore by Proposition
 \ref{prop:celgraph} it gives the representation of $\Cc_K$ which, explicitly,  consisting of:
 
 \begin{itemize}
 \item The stalks $E_x, E_e$ at the vertices and edges of $K$. We write $E_x(\Fc), E_e(\Fc)$ if needed.
 
 \item The generalization maps $\gamma_h: E_x\to E_e$ for any incidence, i.e. half-edge $h$ containing $x$ and
 contained in $e$.

 \end{itemize}

 \begin{prop}
 For the Verdier dual perverse sheaf $\Fc^\bigstar$ we have canonical identifications
 \[
 E_x(\Fc^\bigstar) \simeq E_x(\Fc)^*, \quad E_e(\Fc^\bigstar) \simeq E_e(\Fc)^*.
 \]
 \end{prop}
 
 \noindent{\sl Proof:} Follows from the local statement for a star shaped graph in a disk, Prop. \ref{prop:verdier}.\qed

 \vskip .2cm
 
 So we define
 \[
 \delta_h = (\gamma_h^{\Fc^\bigstar})^*: E_e\to E_x. 
 \]
 
 \begin{thm}\label{thm:main}
 (a) The data $Q=Q(\Fc) = (E_x, E_e, \gamma_h, \delta_h)$ form an object of the category $\Ac_K$. 
 
 (b) If $K$ is a spanning graph for $S$, and $N=\Vert(K)$, then the functor 
 \[
 Q_K: \Perv(S,N) \lra\Ac_K, \quad \Fc\mapsto Q(\Fc)
 \]
 is an equivalence of categories. 
 
 \end{thm}
 
 \noindent {\sl Proof:} (a) The relations ($C1_x$)-($C3_x$) resp. ($C_x$) defining $\Ac_K$, are of local nature,
 so they follow from the local statement (Proposition \ref{prop:rels}) about a star shaped graph in a disk. 
 
 \vskip .2cm
 
 (b) This is obtained by gluing the local results (Theorem \ref{thm:An}). More precisely, perverse sheaves smooth
 outside $N$, form a stack $\Pen$ of categories on $S$. We can assume that $S=\Surf(K)$, so $\Pen$ can be seen
 as a stack on $K$, and $\Perv(S,N) = \Gamma(K, \Pen)$ is the category of global sections of this stack. 
 Similarly, $\Ac_K$ also appears as $\Gamma(K, \Aen)$, where $\Aen$ is the stack of categories on $K$ given
 by $K'\mapsto \Ac_{K'}$ (here $K'$ runs over open subgraphs of $K$). Our functor $Q$ comes from a morphism
 of stacks $\Qen: \Pen\to\Aen$, so it is enough to show that $\Qen$ is an equivalence of stacks. This can
 be verified locally, at the level fo stalks at arbitrary points $x\in K$, where the statement reduces to
 Theorem \ref{thm:An}. \qed
 
 \vskip .2cm
 
 Cf. \cite{KS}, \S 9B for a similar argument. 
 
 \vskip .2cm
 
 \noindent {\bf D. Polarized sheaves.} Let us call a {\it polarized} 
 space an object $E\in \Vect(\kb)$ equipped with a nondegenerate $\kb$-bilinear form
$$
\langle \,\,  , \,  \rangle:\ E\times E\lra \kb,
$$
not necessarily symmetric. 
A linear map
$f:\ E\to E'$ between polarized spaces has  
 two adjoints: the left and the right one, 
$^\top f,  f^\top:\ E'\to E$, defined by
$$
\langle {}^\top{} f(x), y\rangle = \langle x, f(y)\rangle;\ 
\langle y, f^\top(x)\rangle = \langle f(y), x\rangle.
$$
Polarized spaces give rise to several interesting geometric structures motivated by category theory, see \cite{Bon}.  

Let us call a {\it polarized perverse sheaf} over $S$ an object $\Fc\in \Perv(S,N)$ equipped with an isomorphism with its Verdier dual
$B: \Fc\to \Fc^\bigstar$.  This concept can be compared with  that of  \cite{saito};
however we do not require any symmetry of   $B$. 
 Polarized perverse sheaves on $S$ with singularities in $N$ form a category $\PolPerv(S,N)$
 whose morphisms are morphisms of perverse sheaves commuting with the isomorphisms $B$. 

\begin{Defi}
Given a ribbon graph $K$, we define the category $^p\Ac_K$ 
(resp. $\Ac^p_K$) whose objecs are collections $Q=(E_x, E_e, \gamma_h, \delta_h)$ as in $\Ac_K$, together 
with an additional  data of polarizations on
  all spaces $E_e, E_x$, subject to an additional condition: 
 for each $h$,  we have
$\delta_h = {}^\top\gamma_h\ (\text{resp.}\   \delta_h = \gamma_h^\top).$
\end{Defi}

Since by definition the equivalence $Q_K$ from Theorem \ref{thm:main} commutes 
with the duality,  we obtain:

\begin{cor} If $K$ is a spanning graph for $S$ and $N=\Vert(K)$,  then
we have two equivalences
$$
^pQ_K:\ \PolPerv(S,N) \iso \ ^p\Ac_K;\quad  Q_K^p:\ \PolPerv(S,N) \iso \ \Ac_K^p. \qed
$$
\end{cor}

 \appendix
 
 \section{Appendix. Coboundary actions and helices}
 
 Since   the same marked surface $(S,N)$ has many spanning graphs $K\supset N$,  the corresponding categories $\Ac_K$
  are all equivalent to each other, being identified with $\Perv(S,N)$. 
 One can continue the analysis of this paper by constructing  a system of explicit identifications $\Ac_K\to\Ac_{K'}$ for pairs of
 spanning graphs $K, K'\supset N$ connected by ``elementary moves'', in the spirit of \cite{DK}. 
 To keep the paper short, we do not do it here, but  discuss a local aspect of this issue: the action of $\ZZ/n$ on $\Ac_n$. 
 
 \vskip .2cm
 
 Let $G$ be a discrete group which we consider as a category with one object $\pt$. Recall (see, e.g., \cite{De}, \cite{GK})  
 that a {\em category with $G$-action}, or a {\em categorical representation} of $G$ is a lax 2-functor $F: G\to\Cc at$
 from $G$ to the 2-category of categories. 
 
 Explicitly, it consists of the following data (plus the data involving the unit of $G$, see \cite{GK}): 
  
 \begin{itemize}
\item[(0)]  A category $\Cc= F(\pt)$.

\item[(1)]  For each $g\in G$, a functor $g_*=F(g): \Cc \to\Cc$. 

\item[(2)] For each $h,g\in G$, an isomorphism of functors $\alpha_{h,g}= F(h,g): h_* g_* \Rightarrow (hg)_*$.

\item[(3)] It is required that for any three elements $h,g,f\in G$ the square 
$$
\begin{matrix} h_*g_*f_* & \overset{\alpha_{h,g}}\lra & (hg)_*f_*\\
\alpha_{g,f}\downarrow & &\downarrow \alpha_{hg,f}\\ 
h_*(gf)_*& \overset{\alpha_{h,gf}}\lra & (hgf)_* 
\end{matrix}
$$
is commutative, i.e.   
$$
\alpha_{hg,f}\alpha_{h,g} = \alpha_{h,gf}\alpha_{g,f}.
$$
 \end{itemize}
 
\begin{ex}\label{ex:cocycle}  If $F(g) = \Id_\Cc$ for all $g\in G$,  then $F$ is the 
same as a $2$-cocycle $\alpha\in Z^2(G;Z(\Cc))$ where $Z(\Cc)$ is the 
{\it center} of $\Cc :=$ the group of automorphisms of the identity functor $\Id_\Cc$, the action of $G$ on $Z(\Cc)$ being trivial. 
\end{ex}

We say that the action is {\it strict} if $(gf)_* = g_*f_*$,  
and $\alpha_{g,f} = \Id_{g_*f_*}$ 
for all composable $g, f$. In other words, a strict action is simply 
a group homomorphism
$
F: G \to \Aut(\Cc). 
$

\vskip .2cm

All actions of $G$ on a given category $\Cc$ form themselves a category, denoted $\Act(G, \Cc)$. It has a distinguished object $I$,
the trivial action, with all $g_*=\Id_\Cc$ and all $\alpha_{h,g}=\Id$. 
 Given an action $F$ as above, {\it a coboundary structure} on $F$ is an isomorphism $\beta: I\to F$ in $\Act(G,\Cc)$. Explicitly, it consists of:
 
 \begin{itemize}

\item A  collection of natural transformations
$
\beta_g: \Id_\Cc \iso g_*, \ 
$
given for all $g\in G$
such that 

\item  For each  $g, f\in G$ the square 
$$
\begin{matrix} \Id_\Cc& \overset{\beta_f}\lra & f_*\\
\beta_{gf}\downarrow & & \downarrow \beta_g\\
(gf)_*& \overset{\alpha_{g,f}^{-1}}\lra & g_*f_*
\end{matrix}
$$
is commutative, in other words, 
$
\alpha_{g,f}= \beta_{gf}\beta_f^{-1}\beta_g^{-1}. 
$

\end{itemize}

\begin{exas} (a) 
 In the situation of Example \ref{ex:cocycle} a coboundary structure on $F$ is the same as a $1$-cochain
\[
\beta\in C^1(G;Z(\Cc)) = \Hom_{\on{Set}}(G, Z(\Cc)), \quad d\beta = \alpha.
\]
 
 (b) If our action is strict, 
then  a coboundary structure on it is a collection of natural 
transformations $\{\beta_g\}$ as above, such that 
$
\beta_{gf} = \beta_{g}\beta_{f}.
$
\end{exas}

Returning now to the situation of \S \ref {sec:fracspin}, we have a strict action
 of $\ZZ/n$ on $\Ac_n$ such that $k\in\ZZ/n$ acts by
 rotation by $2\pi k/n$.  More precisely, for 
$$
x = (E_0, E_1, \ldots, E_n;\ \gamma_i, \delta_i)\in \Ob(\Ac_n)
$$
 we define
$$
k_*x = (E_0, E_{1+k},\ldots, E_{n+k};  \gamma_{i+k}, \delta_{i+k}),\quad
k\in \ZZ/n
$$
where the indices (except for $E_0$) are understood modulo $n$. 

\begin{prop}  The strict action of $\ZZ$ on $\Ac_n$ induced by the composition
 \[
\ZZ \lra \ZZ/n \lra Aut(\Ac_n)
 \]
  is coboundary.
  
  \end{prop}

  \noindent {\sl Proof:} For $x\in \Ac_n$ as above 
we have an arrow $\beta_1(x):\ x \to 1_* x$
in $\Ac_n$, induced be the fractional monodromies
$$
T_i = \gamma_{i+1}\delta_i: E_i \lra E_{i+1},\ i\in \ZZ/n,
$$ 
which give rise to  natural isomorphisms 
$\beta_1:\ \Id_{\Ac_n} \to 1_*$
(here $1\in \ZZ/n$ is the generator). 
More generally, putting 
$
\beta_k := (\beta_1)^k,\ k\in \ZZ,
$
we get a coboundary structure on the  composed action.  \qed

\vskip .2cm

 Note that in particular the global monodromy
$T = \beta_n$ is a natural transformation 
$
\Id_{\Ac_n} \buildrel\sim\over\to  \Id_{\Ac_n},
$
so it is an element of $Z(\Ac_n)$. 
For an object $x\in\Ac_n$ the sequence
\[
\cdots (-1)_*x, x, 1_*x, \cdots, n_*x, \cdots
\]
can be seen as a decategorified analog of a {\em helix}, see \cite{BP}, 
with the monodromy $T$ playing the role of the Serre functor.

 \let\thefootnote\relax\footnote {
M.K.: Kavli Institute for Physics and Mathematics of the Universe (WPI), 5-1-5 Kashiwanoha, Kashiwa-shi, Chiba, 277-8583, Japan; 
mikhail.kapranov@ipmu.jp

%Email: {mikhail.kapranov@ipmu.jp}

\vskip .2cm

V.S.: Institut de Math\'ematiques de Toulouse, Universit\'e Paul Sabatier, 118 route de Narbonne, 
31062 Toulouse, France; schechtman@math.ups-toulouse.fr 

}

\end{document}